%Paper: q-alg/9507011
%From: "A.Kazarnovski-Krol" <akrol@math.rutgers.edu>
%Date: Mon, 17 Jul 1995 02:50:35 -0400

%this is Ams-Tex

\input amstex.tex
\documentstyle{amsppt}
\input amsppt1.tex
\nologo
\magnification=1200
%\nopagenumbers

  %\parindent=0pt
\NoBlackBoxes

\define\op#1{\operatorname{#1}}

\topmatter
\author   A. Kazarnovski-Krol            \endauthor
\title A generalization of Selberg integral \endtitle
\address{
  Department of Mathematics
  Rutgers University
  New Brunswick, NJ 08854, USA}

\vskip 10mm

\abstract{ We analyze the situation which is related to zonal
spherical functions of type $A_n$
 and obtain  a generalization of Selberg integral.
}
\endtopmatter

\document

\subhead 0.1 Notations \endsubhead

$R$ -set of roots of root system of type $A_n$

$R_{+}$ - set of positive roots

$\alpha_1, \alpha_2, \ldots, \alpha_n$ - simple roots of root system
of type $A_n$

$\delta= {1\over 2} \sum_{\alpha \in R_{+}} \alpha$ -halfsum of
positive roots

$ \Lambda_1, \Lambda_2, \ldots, \Lambda_n$ -fundamental weights, i.e.

$(\alpha_i, \Lambda_j) = \delta_{ij}$ , where $\delta_{ij}$ is
Kronecker's delta

$k$- complex parameter ( 'halfmultiplicity' of a root)
$$\rho =\rho(k)= {k \over 2} \sum_{\alpha \in R_{+}} \alpha$$
Let $\Bbb R ^{n+1}$ be a $(n+1)$-dimensional Euclidean vector space
with inner product $(.,.)$ and with basis $e_1, e_2, \ldots e_{n+1}$.
We realize simple roots $\alpha_1, \alpha_2, \ldots, \alpha_{n}$
as $e_1 -e_2, e_2 -e_3, \ldots, e_n - e_{n+1}$, correspondingly.

$\alpha^{\vee}=\frac{ 2 \alpha}{(\alpha, \alpha)}$
\smallskip

$W=S_{n+1}$ Weyl group of type $A_n$ (group generated by the
orthogonal reflections with respect to hyperlanes perpendicular to
roots $\alpha \in R$)

$\{z_l, \;\;l=1, \ldots, n+1\}$ - arguments

$\{ t_{ij} ,\; i=1, \ldots, j,\;\;
j=1,\ldots, n\}$ - variables of integration

$$\lambda=(\lambda_1,\ldots, \lambda_{n+1})\; |\lambda_1 +\lambda_2
+\ldots +\lambda_{n+1}=0$$
Though this homogeneity condition might be released we prefer to
impose it. Also parameter $\lambda$ is assumed to be generic.

$\phi(\lambda +\rho(k),k,z)$- asymptotic solution with the leading
asymptotic $z^{\lambda +\rho}$, i.e.

$\phi(\lambda+\rho, k, z)= z^{\lambda+\rho}(1+ \ldots)$

$\Delta_w(z)\;, w \in W$ - cycles for asymptotic solutions, cf. ref.[43]

 \head 1. Introduction  \endhead

 We analyeze the situation which is related to zonal spherical
functions of type $A_n$. The situation is also known as
Calogero-Sutherland model.
The zonal spherical functions on symmetric Rimannian spaces were
introduced in ref. [70].
In the case of root system of type $A_n$ it is by proven by
I.Cherednik and
A.Matsuo in refs.
[38,39] that hypergeometric system of
differential equations of Heckman-Opdam cf. ref. [45] is related to
the particular
case of
trigonometric version of
Knizhnik-Zamolodchikov equation in conformal field theory.
In particular, solutions to the hypergeometric system of Heckman-Opdam
can be obtained from the solutions of Knizhnik-Zamoldchikov equations
 by symmetrization procedure.
Solutions to Knizhnik-Zamolodchikov equations  are given by certain
 multidimensional integrals, whose
integrand has the standard part times complicated meromorphic
function cf. refs. [40, 32].
This complicated meromorphic factor becomes even more complicated
(formally) after symmetrization. We would like  like to emphasize that
in this particular case
one
can get rid of this unpleasent meromorphic factor cf. theorem
3.2. below.
See also the refs. [16 ,68 ] , as well as refs. [19,15,69] about W-algebras.
\smallskip

Here is the organization of the paper.
In section two we recall the transformation law for the Heckman-Opdam
hypergeometric functions related to root system cf. ref.[48].
In section 3  using the integral
representation in the case of root system of type  $A_n$ from
ref. [43]
with the help of transformation law we obtain another integral
representation (theorem 3.2)  and calculate the leading coefficient.
In section 4 using the evaluation theorem of Opdam ( theorem 4.4 below)
we obtain generalized Selberg integral (theorem 4.1).
Remarkably, the answer is the same (up to the phase)
for  $(n+1)!$ different contours of integration.

 Recently integrals of this type have drawn much attention because of
applications to conformal field theory, cf. ref. [2].

\smallskip
 \head 2. Transformation law              \endhead

\subhead 2.1. Differential operator of second order \endsubhead

Let $L$ be the following differential operator
$$
L= L(k)=
{\sum_{i=1}^{n+1} (z_i \frac{\partial} {\partial z_i})^2} -
{k\sum_{i < j} {\frac{z_j +z_i}{z_j-z_i}}(z_i \frac{\partial} {\partial z_i}-
z_j \frac{\partial} {\partial z_j})}.
$$
\remark{Remark 2.2 }
Operator $L$ originates in the theory of zonal spherical functions as
the radial part of Laplace-Casimir operator of second order
taken with respect to
Cartan decomposition $G=KAK$ cf. refs.[46,9].
\endremark
\smallskip

\subhead{2.3 Important  property of operator $L$} \endsubhead
$$
\prod (z_i-z_j)^{2 k-1}\prod z_i^{(1- 2k) \frac{n}{2}} \circ
L(k)  \circ \prod (z_i-z_j)^{1-2 k } \prod z_i^{(2 k -1) \frac{n}{2}}=
L(1-k) +(1 -2 k)(\delta, \delta)
$$
cf. Proposition 4.2. ref. [48].
We recall that $\delta$ is half the sum of positive roots:
$\delta= {1\over 2} \sum_{\alpha \in R_{+}} \alpha$.
\smallskip

\subhead {2.4 Asymptotic solutions}
\endsubhead
For generic $\lambda$
operator $L=L(k)$ has $(n+1)!$ eigenfunctions with the eigenvalue $(\lambda,
\lambda) -(\rho,\rho)$ with leading asymptotic $z^{w \lambda +\rho}$,
correspondingly.
Recall that $\rho= {k\over 2} \sum_{\alpha \in R_{+}} \alpha$.
Asymptotic solutions are enumerated by the elements of the Weyl
group $w \in W$.
These solutions satisfy to the whole hypergeometric system of
differential equations of Heckman and Opdam (Macdonald, Sekigushi,
Debiard). Moreover, locally they provide a basis for all the solutions of the
hypergeometric system cf. Corollary 3.11 ref. [45].
We denote asymptotic solution with leading asymptotic $z^{\lambda +\rho}$ by
$\phi(\lambda+\rho(k),k,z)$.
Asymptotic solutions are connected with many interesting parts of
mathematics and physics, cf. refs. [ 23,36 ,37].

\proclaim {Theorem 2.5 }(Transformation law)(Opdam)
$$
\prod (z_i- z_j)^{2 k -1} \prod z_i^{(1 -2 k){n \over 2}}
\phi(\lambda+\rho(k), k,z) =\phi(\lambda +\rho(1-k), 1-k,z)
$$
cf. Corollary 4.4 ref. [48].
\endproclaim
 The proposition is an easy corollary of
the property 2.3. of operator $L$.

\remark{2.6 Note:} The homogeneity of
$$\prod(z_i-z_j)^{1-2k} \prod z_i^{(2 k -1){n \over 2}}$$
equals to zero.
\endremark

\remark{Remark 2.7} The importance of the above simple theorem 2.5 is hardly
possible to overvalue.
\endremark

\head 3. Integral representations \endhead

Consider the following variables: $z_l,\; l=1,...,n+1;\;
t_{i,j},\;i=1,...,j,\;j=1,...,n$ . It is
convenient to organize these variables in the form
of a pattern cf. fig. 1. The idea of such an organization is borrowed
from ref.[8], while variables $t_{ij}$ itself have a nice geometric
origin in elliptic coordinates cf. ref. [10]. Also these variables
appear in Knizhnik-Zamolodchikov approach in ref. [38].

\midinsert
$$
\matrix
z_{1}&&z_{2}&&\ldots&&\ldots&&z_{n+1}\\\\
&t_{1,n}&&t_{2,n}&&\ldots&&t_{n,n}&\\\\
&&\ldots&&\ldots&&\ldots&&\\\\
&&&t_{1,2}&&t_{2,2}&&&\\\\
&&&&t_{1,1}&&&&
\endmatrix
$$
\botcaption{Figure 1}
 Variables organized in a pattern
\endcaption
\endinsert

\smallskip
   In ref. [43] we described contours for integration $\Delta_w= \Delta_w(z)$
which provide asymptotic solutions $\phi(w \lambda +\rho,k,z)$ for the
Heckman-Opdam hypergeometric system of differential equations.
We also obtained a multivalued form
and made a natural convention about the phase of the form over cycle
$\Delta_w$.
We assume that similar convention is made in theorems 3.2 and 4.1 below.

\proclaim{ Theorem 3.1}
 Let $w\in{S_{n+1}} $. Then for generic $\lambda,\; k$ the integral of the
multivalued
form below over cycle $\Delta_w$ gives an asymptotic solution
$\phi(w \lambda+\rho, k,z) $
to the Heckman-Opdam hypergeometric system of differential equations

$$
\align
 &
 \prod_{i=1}^{n+1} z_{i}^{\lambda_1 +{k n\over 2}}
 \prod_{i_1 > i_2}(z_{i_1}-z_{i_2})^{1-2k}
\int_{\Delta_w(z)}
 \prod_{i,i_1} ( z_i - t_{i_1,n} )^{k-1}\\
 &\times \prod_{j=1}^{n-1} \prod_{i,i_1} ( t_{ij} - t_{i_1,j+1} )^{k-1}\\
 &\times \prod_{j=2}^{n} \prod_{i_1>i_2} {(t_{i_1,j}-t_{i_2,j})^{2-2k} }\\
 &\times\prod_{j=1}^n \prod_{i=1}^{j} {t_{ij}^{\lambda_{n-j+2}-
  \lambda_{n-j+1} -k} } \quad
 {dt_{11} dt_{12} dt_{22} \ldots dt_{nn} }=
 a(w) \phi(w \lambda+\rho, k,z)
\endalign
$$

$$
\align
a(w)=
& e^{-2\pi i(\lambda, \delta)} e^{- \pi i (k-1) l(w)}
 (2i)^{\frac{n(n+1)}{2}} \Gamma(k)^{\frac{n(n+1)}2}\\
&\times \prod_{\alpha\in R_{+}}
\frac{\Gamma((-w{\lambda}, {\alpha^\vee}))
 \sin( \pi(-w \lambda, \alpha^\vee))}
{\Gamma((-w{\lambda},{\alpha^\vee)}+k)}
\endalign
$$

cf. theorems 6.1 and 6.3 of ref. [43].
\endproclaim
\smallskip
Integral representation in theorem 3.1 has certain advantage. Namely,
one can identify contour for integration for zonal spherical function
itself cf. ref.[44].
\smallskip

Applying the transformation law to the integral representation of
theorem 3.1.
one obtains the following integral representation.

\proclaim{ Theorem 3.2}
 Let $w\in{S_{n+1}} $. Then for generic $\lambda,\;k$ the integral of
the
multivalued
form below over cycle $\Delta_w$ gives an asymptotic solution
$\phi(w \lambda+\rho, k,z) $
to the Heckman-Opdam hypergeometric system of differential equations:

$$
\align
 \prod_{i=1}^{n+1} z_{i}^{\lambda_1 +{k n\over 2}}
 &\int_{\Delta_w(z)}
 \prod_{i,i_1} (z_i -t_{i_1,n})^{-k}\\
 & \times\prod_{j=1}^{n-1} \prod_{i,i_1} ( t_{ij} - t_{i_1,j+1} )^{-k}\\
 &\times\prod_{j=2}^n  \prod_{i_1 > i_2} (t_{i_1,j}-t_{i_2,j})^{2k}\\
 &\times \prod_{j=1}^n \prod_{i=1}^{j} {t_{ij}^{\lambda_{n-j+2}-
  \lambda_{n-j+1} +k-1} } \quad
 {dt_{11} dt_{12} dt_{22} \ldots dt_{nn} }=
 a(w)\phi(w \lambda+\rho, k,z)
\endalign
$$
where
$$
\align
a(w)=
& e^{-2\pi i(\lambda, \delta)} e^{\pi i k l(w)}
 (2i)^{\frac{n(n+1)}{2}} \Gamma(1-k)^{\frac{n(n+1)}2}\\
&\times \prod_{\alpha\in R_{+}}
\frac{\Gamma((-w{\lambda}, {\alpha^\vee}))
 \sin( \pi(-w \lambda, \alpha^\vee))}
{\Gamma((-w{\lambda},{\alpha^\vee)}-k +1)}
\endalign
$$
\endproclaim
\smallskip

\remark{Remark 3.3} Calculation of leading asymptotic coefficient uses
diagrams cf. ref. [43], section 1,2  and induction. Or can be obtained from
theorem 3.1 simply by replacing $k$ by $1-k$.
\endremark
\smallskip

\remark{Remark 3.4} Compare integral representation in theorem 3.2 with
the integral representation indicated in ref. [38], obtained with the
help of symmetrization of solutions of trigonometric
Knizhnik-Zamolodchikov equation. Also, cf. ref. [40] for the integral
solutions
of trigonometric Knizhnik-Zamolodchikov equations.
\endremark

\smallskip

\head{4. Generalized Selberg integral}\endhead

Here is the main result of the paper.

\proclaim {Theorem 4.1 }( Generalized Selberg integral)
$$
\align
& \int_{\Delta_w(1)}
 \prod_{i_1=1}^n (1 -t_{i_1,n})^{-(n+1)k}\\
 & \times\prod_{j=1}^{n-1} \prod_{i,i_1} ( t_{ij} - t_{i_1,j+1} )^{-k}\\
 &\times\prod_{j=2}^n \prod_{i_1 > i_2} (t_{i_1,j}-t_{i_2,j})^{2k}\\
 &\times\prod_{j=1}^n \prod_{i=1}^{j} {t_{ij}^{\lambda_{n-j+2}-
  \lambda_{n-j+1} +k-1} } \quad
 {dt_{11} dt_{12} dt_{22} \ldots dt_{nn} } \\
& = {(2 \pi i)^{n(n+1)\over 2}} e^{-2 \pi i(\lambda, \delta)}
 e^{\pi i k l(w)}
{\Gamma(1-k)}^{n(n+1)\over2}\\
&\times \prod_{\alpha \in{R_{+}}}
 \frac{1}{\Gamma((w \lambda, \alpha^{\vee})-k+1)
\Gamma((-w \lambda, \alpha^{\vee})-k+1)}\\
&\times \prod_{\alpha \in{R_{+}}} \frac{\Gamma((-\rho, \alpha^{\vee})-k+1)}
     {\Gamma((-\rho, \alpha^{\vee})+1)}
\endalign
$$
\endproclaim

\remark{Remark 4.2} Remarkably  the above constant does  depend on
$w \in S_{n+1}$ only in the phase factor $e^{ \pi i k l(w)}$.
Note also, that in ref. [43] we made a natural convention about the
phase of the integrand $\omega_w$ over $\Delta_w$.
\endremark
\smallskip
\remark{Remark 4.3}
The generalized Selberg integral can be conveniently rewritten
as follows. Let's assign to each variable $t_{ij}$ a simple root
$\alpha(t_{ij})$ by
the rule:

$$\alpha (t_{ij}) = \alpha (j)= \alpha_{n +1-j}= e_{n+1-j} - e_{n+2-j}$$
Note that to each variable of the same row we assign the same simple
root. This assignment looks different from [38] only  because we use
different indexation of variables of integration $t_{ij}$.

Let also $\Lambda_1$ be the first fundamental weight.

Then one can rewrite the Selberg integral from theorem 4.1 as follows:
$$
\align
& \int
\prod t_{ij}^{(\lambda-\rho, -\alpha(j))}\\
&\times \prod (1- t_{ij})^{k((n+1) \Lambda_1, -\alpha(j))}\\
& \times \prod (t_{ij} - t_{i'j'})^{k(-\alpha(j), -\alpha(j'))}
\frac {dt_{11}}{t_{11}} \frac {dt_{12}}{t_{12}} \ldots  \frac
{dt_{nn}}{t_{nn}}
\endalign
$$
\endremark

The theorem is an immediate application of the theorem 3.2 and the
following theorem.
\smallskip
\proclaim{Theorem 4.4} Evaluation theorem (Opdam).
$$
\phi(w \lambda+\rho(k),k,1)=
\lim_{z \rightarrow 1} {\phi(w \lambda+\rho(k),k,z)}=
{\frac
{\prod_{\alpha \in R_{+}}
      {\frac {\Gamma( (w \lambda, \alpha^{\vee})+1)}
      {\Gamma( (w \lambda, \alpha^{\vee})- k +1)}}}
{ \prod_{\alpha \in R_{+}}
 {\frac
 {\Gamma( -(\rho , \alpha^{\vee})+1)}
      {\Gamma( -( \rho, \alpha^{\vee})- k +1)}}}}
$$

cf. theorem 6.3 [48].
\endproclaim
Recall that we restrict ourselves to the case of root system of type
$A_n$.
\smallskip

\remark{Concluding remarks   }
In this paper we
obtained a generalization of Selberg integral.
Integrals of this type
are important for the conformal field theory cf. ref. [2].
The Selberg integral considered in this paper can serve as an
example of such integrals.

\endremark

\subhead Acknoledgements \endsubhead
I am greatful  to  I. Gelfand and  S. Lukyanov for  stimulating
discussions.

\end